\documentclass[12pt,titlepage,a4paper]{article}
\usepackage{amssymb,latexsym,amsfonts,amsmath,amsthm} 
\usepackage{verbatim}

{\makeatletter \@addtoreset{equation}{section}}
\newcommand{\blackbox}{\hfill \rule{2mm}{2mm}}

\newcommand{\Ker}{\mathop{\mathrm{Ker}}\nolimits}
\newcommand{\IIm}{\mathop{\mathrm{Im}}\nolimits}

\newcommand{\mcX}{\mathop{\mathcal{X} {} }\nolimits}

\newcommand{\mcU}{\mathop{\mathcal{U} {} }     \nolimits}
\newcommand{\mcV}{\mathop{\mathcal{V} {} }\nolimits}

\newcommand{\h}{ \mathop{ \mathrm{h} {} }\nolimits }
\newcommand{\e}{ \mathop{ \mathrm{e} {} }\nolimits }
\newcommand{\g}{ \mathop{ \mathrm{g} {} }\nolimits }

\newtheorem{theorem}{Theorem}[section]
\newtheorem{lemma}[theorem]{Lemma}

\newcommand{\Proof}{ \noindent{\bf Proof:~} }

\begin{document}

\title{Regular   submodules  of torsion
modules over a discrete valuation domain } 
\author{Pudji Astuti
\thanks{The research of the first author was
supported by Deutscher Akademischer Austauschdienst under
Award No. A/98/25636.}
\\
Departemen Matematika\\
Institut Teknologi Bandung\\
Bandung 40132\\
Indonesia
\and 
Harald K. Wimmer\\
Mathematisches Institut\\
Universit\"at W\"urzburg\\
D-97074 W\"urzburg\\
Germany}
%


\maketitle

\begin{abstract}

A submodule $W$ of a p-primary module $M$ of bounded order 
is known to be regular if $W$ and $M$ have simultaneous bases.
In this paper we derive necessary and sufficient conditions for
 regularity of a submodule. 

\vspace{2cm}
\vspace{.5cm}


\end{abstract}

\section{Introduction} \label{sec.1}
Let $R$ be a discrete valuation domain
with maximal ideal $Rp$, and let $M$ be a torsion module over $R$
 and $W$  be a submodule of $M$. 
The submodule $W$ is called \emph{regular } \cite[p.65]{Kap}, 
\cite[p.102]{Vil} if
\begin{equation} \label{e.kapl}
 p^n W \cap p^{n+r}M = p^n(W \cap p^r M)
\end{equation}
holds for all $n \geq 0$,  $r \geq 0$. The 
regularity condition \eqref{e.kapl} was introduced by
Vilenkin \cite{Vil} in his study of 
decompositions of topological p-groups. Kaplanski 
\cite{Kap} showed that for a 
module $M$ of bounded order  \eqref{e.kapl} is necessary and 
suffifient for the existence of simultaneous bases of
$W$ and $M$.  
In this paper we 
shall identify two conditions which are equivalent to
\eqref{e.kapl}.
One is related to a theorem of Baer 
\cite[p.4]{FuII} on the decomposition of elements in
abelian p-groups, the other one was introduced by 
Ferrer, F. Puerta and X. Puerta
\cite{FPP} to characterize marked invariant subspaces
of a linear operator.

Notation and definitions:
 The letters $\mcU, \, \mcV , \, 
\mcX , \dots , $ will always denote subsets of $M$.
Let 
$ \langle \mcX \rangle $ be the submodule spanned by
$\mcX  $.
 We shall
use the letters $u, v, x, \dots ,$ for elements of the
module $M$, and $\alpha, \beta, \mu , \dots ,$ will be 
elements of the ring $R$.
Using   the terminology for abelian $p$-groups in 
\cite[p.4]{FuI} 
we say that  $x \in M$ has  \emph{exponent} $k$,
 and we 
write $\e (x) = k$, if  $k$ is the smallest nonnegative 
integer 
 such that $p^{k}x = 0$.
An  element $x \in M$ is said to have (finite) 
\emph{height} $s$ if $x \in p^s M$ and $x \notin p^{s+1} M$,
and  $x$ has \emph{infinite height}, if 
 $x \in p^s M$ for all $s \geq 0$. 
We write  $\h (x)$ for the height of  $x$. 
If $x \in W$ then $\h _W (x)$ will denote the
height of $x$ with respect to $W$.
Note that $\e (0) = 0$ and $\mathrm{h} (0) = \infty$.
 Let $R^{\ast}$ be the group
of units of $R$. 
If $\alpha \in R$ is nonzero and  $\alpha = p^s \gamma$, 
$\gamma \in R^{\ast}$, then we set $\h (\alpha)  = s $.
We put $\h (\alpha)  = \infty$ if 
$\alpha = 0$. 
We call  $x \in M$  an $(s,k;s_1)$-element if $x \neq 0$ and
\[
\h(x) = s, \, \e(x) = k, \, \h(p^{k-1} x) = (k-1) + s_1.
\]
In accordance with a definition of Baer \cite {Bae1} 
we say that an element $x$  is \emph{regular }
if $\h (x) = \infty$ or if $\h (x)$ is finite and 
 \begin{equation}  
\label{e.preg}
\h (p^j x) = j + \h (x), \,\, 
j = 1,\dots , \e(x)-1.
\end{equation}
 The two concepts of regularity introduced above are 
consistent. We shall see in Lemma \ref{la.regelemodule}
 that a finite height
element $x \in M$ is regular  if and only if $\langle x \rangle$
is a regular  submodule of $M$.

For $s \geq 0,\, k\geq 0$ we  define the submodules 
    $ M[p^k] = \{ x \in M \, | \, p^k x = 0 \} $
and 
\begin{equation} \label{eq.msk}
        M^s _k = p^s M \, \cap \, M[p^k].
\end{equation}
Then
\[
        M^s _k = \{ x \in M \, | \,  \e (x) \leq k , \,
 \h (x) \geq s     \}.
\]
In particular $M^s _0 = 0$.


Our main result will be the following.
\begin{theorem} \label{thm.m1}
Let $M$ be a torsion module over a discrete valuation
domain
and let $W$ be a  submodule of  $M$.
The following conditions are equivalent.
\begin{description}
\item[(K)]  $W$ is regular , i.e.
if  $n \geq 0, \, r \geq 0$ then
\begin{equation} \label{e.kap1cap}
p^n W \, \cap \, p^{n+r} M = p^n (W \, \cap \, p^r M).
\end{equation}
\item[(B)]
If $x \in W$ is nonzero
then $x$ can be decomposed 
 as 
\begin{equation} \label{eq.full}
x = 
   y_{k_1}^{s _{1} } + \cdots + y_{k_m}^{s _{m} }
\end{equation}
such that 
\[
 y_{k_i}^{s_i} \in W\,\, is \,\, regular ,\,\, 
i = 1, \dots ,m,
\]
and    
\[
   \h ( y_{k_i}^{s_{i} }) = s_{i}, \,\, 
  \e( y_{k_i}^{s_{i} }) = k_{i},
\]
and   
\begin{equation} \label{eq.eee} 
 k_1 > \cdots > k_m > 0 \,\,
and \,\, 
             s _{1} >  \cdots > s_{m}.
\end{equation}
\item[(FPP)] 
If $s \geq 0, \, k \geq 1$, then 
\begin{equation} \label{eq.ffpp}
      (W \, \cap M^{s+1}_k ) + (W \, \cap M^{s}_{k-1} )
  =
     W \, \cap (  M^{s+1}_k +  M^{s}_{k-1} ).
\end{equation}
\end{description}
\end{theorem}

By a result of Baer  \cite[p.4, Lemma 65.4]{FuII} 
 condition (B) is satisfied for
$W= M$. Hence (B) singles out those submodules $W$ where
each element $x \in W$ allows a decomposition \eqref{eq.full}
such that the summands  $y_{k_i}^{s_i}$ can be chosen from
$W$ itself. With regard to 
condition (FPP) we observe that the
inclusion
\begin{equation} \label{eq.nclu}
 (W \, \cap M^{s+1}_k ) + (W \, \cap M^{s}_{k-1} )
  \subseteq
     W \, \cap (  M^{s+1}_k +  M^{s}_{k-1} )
\end{equation}
holds for all submodules $W$.
 
The proof of the theorem will be split into two parts. 
In Section \ref{sec.kb} we show that (B) and (K) are equivalent
and in Section \ref{sec.fppb} we prove the equivalence of 
(B) and (FPP).

\section{Decomposition of elements}  \label{sec.hequb}
We introduce a condition which will be the link between
   (B) and (K) on  one hand and between
(B) and (FPP) on the other.
%
For  a submodule $W$ we define condition
 (H)  as follows. \smallskip \\
{\bf{(H)}} If $x \in W $ is an $(s,k;s_1)$-element 
 then $x$ can be decomposed 
 as
\begin{equation} \label{eq.yzdec}
   x = y^{s_1}_k + z, \,\, y^{s_1}_k \in W, \, 
 z \in W,
\end{equation}
such that
\begin{equation} \label{eq.yhe}
\h( y^{s_1}_k ) = s_1,
 \,\,
 \e(y^{s_1}_k ) = k, \,\, \mbox{and} \,\, \h(z) = s, \,
 \e(z) < k.
 \end{equation}

\medskip
The following technical 
lemma will be useful in several instances. 
It  implies 
 that the element $y^{s_1}_k$ in 
\eqref{eq.yzdec} is regular .

\begin{lemma} \label{la.bound}
Let $x \in M$  be an $(s,k;s_1)$-element. 
Assume 
\begin{equation} \label{eq.sumz}
 x = y + z, \,\, 
z \in M^s _{k-1}.
\end{equation}
Then $y \neq 0$,   $\e(y) = k$, and
\begin{equation} \label{eq.bndd}
s \leq \h(y) \leq s_1.
\end{equation}
 The element
$y$ is regular if and only if $\h(y) = s_1$.
If $x$ is regular then \eqref{eq.sumz} implies $\h(y) = s$.
\end{lemma}

\Proof 
From \eqref{eq.sumz} follows $ p^{k-1} y  = p^{k-1}x
\neq 0$,  and
 $\e(y) = k$. 
Therefore
\begin{equation} \label{eq.lidas}
(k-1) + \h(y) \leq \h( p^{k-1} y ) = \h( p^{k-1} x ) =
(k-1) + s_1,
\end{equation}
which 
yields $\h (y) \leq s_1$. 
It is obvious from
\eqref{eq.lidas} that we have $\h(y) = s_1$
if and only if 
\[
         \h( p^{k-1} y) = (k -1) + \h(y),
\]
i.e., if and only if $y$ is regular. 
If $x$ is regular then $s_1= s$ and \eqref{eq.bndd}
yields  $\h (y)=  s$.
\blackbox

\begin{lemma} \label{la.BeqH}
For a submodule $W$ the conditions $\mathrm{(B)}$ and
$\mathrm{(H)}$  are equivalent.
\end{lemma}

\Proof There is nothing to  prove if $x$ is 
regular. Thus, in the following we assume that 
$x$ is a non-regular element of $ W$ 
with $\h (x) = s$ and $ \e(x) = k$.
In that case we have
$k > 1$, $ s_1 > s$, and   
$\h (p^{k-1} x) =  (k-1) + s_1$.\\
 {\bf{(B)}} $\Rightarrow$  { \bf{(H)}}
Let $x$ be given as in \eqref{eq.full}, 
with $m \geq 2$. 
Put 
$z = 
  y_{k_2}^{s _{2} } + \cdots + y_{k_m}^{s _{m} }$.
Then 
\eqref{eq.eee} implies $\e(z) \leq k_2 < k$
and 
$\h(z) = s_m = s$.
Hence the  decomposition 
$ x =   y_{k_1}^{s _{1} } + z $ is of type {(H) }.\\
 {\bf{(H)}} $\Rightarrow$ {\bf{(B)}}
Let $x$ be an $(s,k;s_1)$-element of $W$ and assume that
$x$ is decomposed according to (H) as 
\begin{equation} \label{eq.immwd} 
x = y^{s_1}_k + z
\end{equation} 
such that \eqref{eq.yhe} holds.
 We know from Lemma \ref{la.bound}
that $y_{k}^{s _{1} }$ is regular. 
Consider $x$ with $s_1 > s, \, k > 1$. Assume as an
induction hypothesis that condition (H) ensures a 
decomposition of type (B) for all 
 $w \in W$ with $\e(w) < k$.
Thus we
have 
\[
        z = z _{l_2} ^{t_2} + \cdots 
       + 
     z _{l_m} ^{t_m}, \,\, m \geq 2,
\]
with properties in accordance with (B).
Thus $\h ( p^{ l_{2} - 1 } z ) = (l_{2} - 1) + t_2$, 
$t_2 \geq s$, 
and 
     $t_2 > \cdots > t_m = s = \h (z)$, 
and
 $k > \e (z) = l_2 > \cdots > l_m > 0$. 
If 
     $s_1 >  t_2$ 
then  we  already have the desired decomposition.
Now suppose 
$ t_2 \geq   s_1$. Let $j$ be such that 
\begin{equation} \label{eq.j}
     t_2 > \cdots t_j \geq s_1 > t_{j+1}.
\end{equation}
Note that  $t_m \geq s_1$ can not occur because of 
$t_m = s$ and $s_1 > s$.
Set
  \[
v = y_k ^{s_1} + ( z _{l_2} ^{t_2} + \cdots 
       + 
     z _{l_j} ^{t_j}).
\]
Then
 $k > l_2$ yields $\e (v) = k$. Since $ y_k ^{s_1}$
is regular we see that 
                          $p ^{k-1} v = p ^{k-1} y_k ^{s_1}$
implies 
$(k-1)+ s_1 = \h (p ^{k-1} v )$. Hence $ \h (v) \leq s_1$.
On the other hand it follows from 
\eqref{eq.j} that
$\h (v) \geq  s_1$. Therefore $\h (v) =  s_1$,
and $v$ is regular. If we 
rewrite \eqref{eq.immwd} in the form
\[
  x = v +  z^{ t_{j+1} } _{ l_{j+1} } + \cdots +
  z^{t_m}_{l_m},
\]
then we have a decomposition with 
  $\h(v) = s_1$ and $  s_1 >  t_{j+1} 
> \cdots >  t _m =s$
and 
  $\e (v) = k > l_{j+1} > \cdots > l_m > 0$.
 \hfill \blackbox

\bigskip

It is not difficult to check that the
 following observation characterizes the numbers $m$,
$k_i$ and
$s _i$ in \eqref{eq.full}.
For a nonzero element 
 $x \in M$
with $\e(x) =k$ define
$\g(x) = \h(x) + \e(x)$.

\begin{lemma} \label{la.unibaer}
Let $x \in M$ be decomposed as
\begin{equation} \label{eq.fl}
x = 
   y_{k_1}^{s_{1} } + \cdots + y_{k_m}^{s_{m} }
\end{equation}
such that 
\[
  \h ( y_{k_i}^{s_{i} }) = s_{i},\,\,
\e( y_{k_i}^{s_{i} }) = k_i, \,\,
\mbox{and} \,\, 
y_{k_i}^{s_{i} } \,\,
\mbox{is regular },\,\, i = 1, \dots ,m.
\]  
and 
\[
k_1 > \cdots > k_m > 0\,\, and \,\,
s_{1} >  \cdots > s_{m}. 
\]
  Set 
   $K = \{k_1, \dots , k_m \}$.
Then $j \in \{1, \dots , k-1\}$,  
is in $K$ 
if and only if
$g(p^jx) > g(p^{j-1}x)$.
Moreover
\[
\h( p^{k_j - 1}x ) = (k_j - 1) + s_j, \,\, j=1,2, \dots m.
\]
In particular, we have $\e(x) = k_1$ and $\h(x) = s_m$.
\end{lemma}


\bigskip   
\section{Equivalence of (K) and (B)} \label{sec.kb}

Condition (K) can be reformulated in a more convenient 
form.

\begin{lemma} \label{la.Kequiv}
We have 
\begin{equation} \label{eq.kapkap}
p^n W \cap p^{n+r}M = p^n(W \cap p^r M), \, n \geq 0, r \geq 0,
\end{equation}
  if and only if
 for each 
$w \in W$ with 
$\h(p^n w ) = n +r $ 
there exists an element 
                          $\tilde{w} \in W$
such that 
\begin{equation} \label{eq.dochnr}
p^n w = p^n \tilde{w} \,\,\,  \mbox{and} \,\,\, \h (\tilde{w}) = r.
\end{equation}
\end{lemma}
\Proof    Obviously \eqref{eq.kapkap}
is equivalent to
 \begin{equation} \label{e.inclK}
p^n W \, \cap \, p^{n+r} M \subseteq p^n (W \, \cap \, p^r M),
\, n \geq 0, \, r \geq 0.
  \end{equation}
 Now \eqref{e.inclK} holds 
if and only if 
\[ 
x \in p^n W, \,\, x \in p^{n+r}M \,\, \mbox{and} \,\,
  x \notin 
            p^{n+r+1}M
\]
imply $ x \in p^n(W \cap p^rM)$. That implication
means the following. If $x = p^n w$ and $ w \in W$ 
and $\h(x) = n+r$, then 
$x = p^n \tilde{w}$ for some $\tilde{w} \in W$ with 
$\h(\tilde{w} ) \geq r$.
Because of $\h( p^n \tilde{w} ) = n +r$ the inequality
$\h(\tilde{w} ) \geq r$ is equivalent to 
$\h(\tilde{w} ) = r$.  \phantom{00} \hfill \blackbox

\begin{lemma} \label{la.regelemodule}
Let $x$ be an element of finite height with $\e(x) = k$.
Then $x$ is regular if and only if the submodule 
$\langle x \rangle$ is regular , i.e.
\begin{equation} \label{eq.xregm}
p^r \langle x \rangle \cap  p^{n+r}M = p^r\bigl( \langle x \rangle
\cap p^nM \bigr),
\, n \geq 0, \, r \geq 0.
\end{equation}
\end{lemma}

\Proof Assume \eqref{eq.xregm}. We want to show that 
$\h(p^{k-1}x) = (k-1) + s_1$ implies $s_1 = \h(x)$.
According to Lemma \ref{la.Kequiv} there exists an element
$\tilde{x} \in  \langle x \rangle$ with properties corresponding
to \eqref{eq.dochnr}, i.e. $\tilde{x} = \gamma p^t x, \,
\gamma \in R^{\ast}$, and $p^{k-1}x = p^{k-1}(\gamma p^t x)$ and
$\h(p^t x) = s_1$. Then we have 
 $t=0$, and $\h(x) = s_1$. 
It is easy to check that
that \eqref{eq.xregm} holds if $x$ is regular . 
\blackbox

\bigskip
{\bf{Proof of Theorem \ref{thm.m1}, Part I: (B) 
$\Leftrightarrow$ (K)} }\\
{\bf{(B)}} $\Rightarrow$  {\bf{(K)}} 
 We want to show that condition (B)
 implies  (K) 
 in the equi\-valent form of Lemma \ref{la.Kequiv}. Let 
$w \in W$ be such that
$\h(p^n w) = n+r$, and  $h(w) = s$, $\e(w) = k_1$.
Then $s \leq r$ and $k_1 > n$.
Hence 
(B) yields a decomposition
\[
   w = y_{k_1} ^{s_1} + \cdots +  y_{k_m} ^{s_m}
\]
where  the elements $y_i ^{s_i} \in W$ are regular ,
 $h(y_i ^{s_i})
= s_i$,  
 and   
\[
          s_1 > \cdots > s_m = s = \h(w)
\]
   and $
          \e(w)=  k_1 > \cdots > k_m  > 0$.
Let $t$ be such that $k_t > n \geq  k_{t+1}$.
Then
\[
 n + r = \h(p^n w) = \h( p^n   y_{k_1} ^{s_1} + \cdots
 +
 p^n  y_{k_t} ^{s_t} ),
\]
and 
$
   \h(p^n w) = \h( p^n  y_{k_t} ^{s_t} ) =  n + s_t$.
Hence $s_t = r$.
Set 
$ 
 \tilde{w} =   y_{k_1} ^{s_1} + \cdots +
y_{k_t} ^{s_t}.
$
Then $\tilde{w} \in W$ and $\h(\tilde{w}) = r$ and 
$p^n w = p^n \tilde{w}$.
\medskip \\
{\bf{(K)}} $\Rightarrow$ {\bf{(B)}}
Because of Lemma \ref{la.BeqH} it suffices to show that
(K) implies (H). Let $x \in W$ be an $(s,k;s_1)$-element.
Set $w = p^{k-1}x$. Then (K), resp. Lemma~\ref{la.Kequiv},
imply that there exists an $\tilde{x} \in W$ such that 
\begin{equation} \label{eq.xtil}
    p^{k-1}x =  p^{k-1} \tilde{x}
\end{equation}
and $\h( \tilde{x}) = s_1$. From 
\eqref{eq.xtil}
follows $\e( \tilde{x}) = k$ and  
$\h(p^{k-1}\tilde{x}) = (k-1) + s_1$. Now
set $z = x - \tilde{x}$. Then \eqref{eq.xtil} 
yields $\e(z) < k$.
Hence
 $x = \tilde{x} + z$ is a   decomposition of  type 
 (H).
\blackbox

\bigskip
As (K) holds for $W = M$  we can 
write each nonzero element $x$ of $M$ according to
(H)  in the 
form \eqref{eq.yzdec}.
Similarly we can decompose $x$ according to (B) as a sum
of the form \eqref{eq.full}.
 In that case we recover the
result of Baer \cite[p.4, Lemma 65.4]{FuII} 
 mentioned in    Section \ref{sec.1}.

\section{Equivalence of (B) and (FPP)} \label{sec.fppb}

In \cite{FPP} J. Ferrer, and F. and X. Puerta studied 
marked invariant subspaces of 
 an endomorphism $A$ of $\mathbb{C}^n$. Their investigation
is based on 
subspaces of the form 
 $ \IIm (\lambda I - A)^s \cap \Ker  (\lambda I - A)^k$.
Thus  the submodules $M^s_k$ in \eqref{eq.msk}
 are a generalization of those
subspaces.
The next lemma is adapted from \cite{FPP}. It 
 characterizes regular elements in terms
of  $ M^s_k $. Note that $ M^s_k \subseteq M^{s_1}_{k_1}$
if $s_1 \leq s$ and $k \leq k_1$. Hence 
$  M^{s+1}_k + M^s_{k-1} \subseteq   M^s_k$.

\begin{lemma} \label{la.hregul}
An element  $x \in M$ satisfies 
\begin{equation} \label{e.xms}
 x \in M^s _k \,\, \mbox{and} \,\, x 
          \notin  M^{s+1} _k + M^s _{k-1}
\end{equation}
if and only if
 \begin{equation} \label{e.skhreg}
 x \,\, {\mbox{is regular }} \,\; and \,\,
\h (x) = s \,\; \mbox{and} \,\, \e (x) = k .
\end{equation} 
\end{lemma}

\Proof ``$\Rightarrow$'' Assume that $x$ satisfies
 \eqref{e.xms}. 
Recall that $x \in  M^s _k$ if and only if both $\h (s) \geq s$
and $ \e (x) \leq k$. Hence 
\begin{equation} \label{eq.nnttin}
   x  \notin  M^{s+1} _k + M^s _{k-1}
\end{equation}
 implies
 $\h (x) = s$
and $ \e (x) = k $.  
Assume 
$\h (p^{k-1} x) = (k-1) + s_1$.
If we decompose $x$ according to (H)  then
$x = y_k^{s_1} + z \in  M^{s_1} _k + M^s _{k-1}$.
Hence \eqref{eq.nnttin} implies
$s_1 = s$, and $x$ is regular .\\
``$\Leftarrow$''
Consider an element $x$ with properties \eqref{e.skhreg}.
Then $x \in  M^{s} _{k}$, 
$x \notin  M^{s+1} _k$ and  $x \notin  M^{s} _{k-1}$.
If 
\[
   x = y + z, \,\, y \neq 0, \,\, z \in  M^s _{k-1},
\]
and $x$ is regular , then it follows from Lemma 
\ref{la.bound} that $\h(y) = s$. Hence we have 
$y \notin  M^{s+1} _k $ and
$ x  \notin  M^{s+1} _k + M^s _{k-1}$.
 \blackbox

\bigskip

{ \bf{Proof of Theorem \ref{thm.m1},
Part II: (H) 
$\Leftrightarrow$  (FPP)}   }\\
{\bf{(H)}} $\Rightarrow$  {\bf{(FPP)}} 
Because of the inclusion \eqref{eq.nclu} 
 the identity \eqref{eq.ffpp}
in (FPP) is equivalent to
\begin{equation} \label{eq.inclu}
 W \cap (M^{s+1} _{k} +  M^{s} _{k-1})  \subseteq
(W \cap  M^{s+1} _{k})  + ( W \cap   M^{s} _{k-1}).
\end{equation}
 We want to show that condition {(H) }
 implies \eqref{eq.inclu} for all $s \geq0$, $k \leq 1$.
Take an element 
\begin{equation} \label{eq.nusspolt}
x \in  W \cap (M^{s+1} _{k} +  M^{s} _{k-1}).
\end{equation}
Then $x \in M^s_k$ and therefore 
$\h(x) \geq s$ and $\e(x) \leq k$. 
To prove that 
\begin{equation} \label{eq.passt}
x \in (W \cap  M^{s+1} _{k})  + ( W \cap   M^{s} _{k-1})
\end{equation} 
we consider three cases. First, let
$\h(x) \geq s+1$
then $x \in W \cap M^{s+1} _{k}$ and \eqref{eq.passt}
is obvious. Secondly, let $\e(x) \leq k-1$. In that case
$x \in W  \cap   M^{s} _{k-1}$.
Now assume $\h(x) = s$ and $\e(x) = k$.
By Lemma \ref{la.hregul} it follows from 
\eqref{eq.nusspolt} that $x$ is not regular. Hence
$\h(p^{k-1}x) = (k-1) + s_1$ and $s_1 > s$. According to
(H)  we have 
$x = y^{s_1}_k + z$ with $y^{s_1}_k \in W \cap  M^{s_1} _{k}$
and $z \in W \cap   M^{s} _{k-1}$, which yields 
\eqref{eq.passt}. \smallskip
\\
 {\bf{(FPP)}} $\Rightarrow$ {\bf{(H)}} Let $x$ be an 
$(s,k;s_1)$-element. 
 If $s_1 = s$ then $x$ is
regular and we have \eqref{eq.yzdec} with $z=0$. Suppose
now that $x$ is not regular , i.e. $s_1 \geq  s+1$. 
Then Lemma \ref{la.hregul} implies 
$x \in  W \cap (M^{s+1} _{k} +  M^{s} _{k-1})$.
From  (FPP) we obtain
\begin{equation} \label{eq.yyyz} 
x = y + z, \,\, y \in  W \cap  M^{s+1} _{k},\,\,
 z \in W \cap   M^{s} _{k-1}.
\end{equation}
Then $y \neq 0$,  $\e(y) = k$ and $\h(y) \geq s+1$.
Let $y$ in \eqref{eq.yyyz} be such that $\h(y)$ is maximal.
We shall see that such a choice of $y$ implies 
$\h(y) = s_1$, and in that case
 \eqref{eq.yyyz} is a decomposition
of type (H). Now suppose that $\h(y) = \tilde{s} < s_1$.
Then, by Lemma~\ref{la.bound}, the element $y \in W$ is not
regular.
Applying Lemma \ref{la.hregul}
to $ y \in  W \cap  M_k^{ \tilde{s} }$
we obtain 
$ y  \in W \cap \bigl(  M_k^{ \tilde{s} +1 }
+
 M^{ \tilde{s} } _{k-1}   \bigr)$.
Thus (FPP) yields 
 \[
   y =  \tilde{y} + z_2, \,\, \tilde{y} 
\in W \cap  M_k^{ \tilde{s} +1 }, \,\, \tilde{y} \neq 0, \,\,
z_2 \in  W \cap   M^{ \tilde{s} } _{k-1}.
\]
Hence 
$
  x = \tilde{y} + (z + z_2),
$
and we have 
  another decomposition of the form
 \eqref{eq.yyyz}, but now with
$h( \tilde{y}) >  \tilde{s}$,  which contradics the maximality
of $\tilde{s}$.
\phantom{000000}\blackbox

\medskip
\textbf{Acknowledgement.} We 
would like to thank  O. Mutz\-bauer  
for a valuable comment.


\begin{thebibliography}{99}


\bibitem{Bae1} Baer, R.: Types of elements
 and the characteristic
subgroups of abelian groups. Proc. London
 Math. Soc. 39(1935), 481--514.



\bibitem{FPP} Ferrer, J.;  Puerta, F.;  Puerta, X.:
 Geometric
characterization and classification of marked subspaces. 
Linear Algebra Appl. 235(1996), 15--34.


\bibitem{FuI} Fuchs,  L.: Infinite Abelian Groups, Vol. I. 
Academic Press, New York 1973.
\bibitem{FuII} Fuchs,  L.:
 Infinite Abelian Groups, Vol. II. 
Academic Press, New York 1973.


\bibitem{Kap} 
 Kaplanski,  I.: 
Infinite Abelian Groups. University of Michigan
Press, Ann Arbor 1954. 

\bibitem{Vil} 
Vilenkin,  N. Ya.: 
 Direct decompositions of topological groups, I.
(in Russian). Mat. Sbornik N. S.
 19(1946), 85--154.



\end{thebibliography}
\end{document}